\documentclass{amsart}
\usepackage{amsfonts, amsthm, amssymb, amsmath, stmaryrd}
\usepackage{mathrsfs,array}
\usepackage{eucal,color,times,enumerate,accents}
\usepackage[all]{xy}
\usepackage{url}
\usepackage{enumitem}
\usepackage{tikz}
\usepackage{tikz-cd}
\usepackage[pdftex,bookmarksnumbered,bookmarksopen]{hyperref}
\usepackage{stackrel}
\hypersetup{colorlinks, linkcolor=blue}

\input xy
\xyoption{all}
\newcommand{\PP}{\mathbb{P}}
\newcommand{\Zl}{\mathbb{Z}_{\ell}}
\newcommand{\Ql}{\mathbb{Q}_{\ell}}

\newcommand{\Gm}{\mathbb{G}_m}

\newcommand{\F}{\mathbb{F}}
\newcommand{\Q}{\mathbb{Q}}
\newcommand{\Z}{\mathbb{Z}}
\newcommand{\N}{\mathbb{N}}
\newcommand{\etale}{\'etale }

\newcommand{\calC}{\mathcal{C}}

\newcommand{\Ker}{\mathrm{Ker}}
\newcommand{\Image}{\mathrm{Im}}

\usepackage{calligra}
\makeatletter
\def\@splitop#1#2\@nil{$\mathscr{#1}\!\!\!$\calligra#2\;\!}
\newcommand*\DeclareCursiveOperator[2]{%
  \newcommand#1{\mathop{\mbox{\@splitop#2\@nil}}\nolimits}}
\makeatother
\DeclareCursiveOperator{\EXT}{Ext}
\DeclareCursiveOperator{\HOM}{Hom}

\begin{document}
\bibliographystyle{alpha}
\newtheorem{theorem}{Theorem}[subsection]
\newtheorem{fact}[theorem]{Fact}

\newtheorem{proposition}[theorem]{Proposition}
\newtheorem{lemma}[theorem]{Lemma}
\newtheorem{corollary}[theorem]{Corollary}
\newtheorem{claim}[theorem]{Claim}
\newtheorem{definition}[theorem]{Definition}

\theoremstyle{definition}
\newtheorem{question}[theorem]{Question}
\newtheorem{conjecture}[theorem]{Conjecture}
\newtheorem{counterexample}[theorem]{Counterexample}

\newtheorem{answer}[theorem]{Answer}
\newtheorem{remark}[theorem]{Remark}
\newtheorem{example}[theorem]{Example}
\newtheorem{warning}[theorem]{Warning}
\newtheorem{notation}[theorem]{Notation}
\newtheorem{construction}[theorem]{Construction}

\newcommand{\adjunction}[4]{\xymatrix@1{#1{\ } \ar@<-0.3ex>[r]_{ {\scriptstyle #2}} & {\ } #3 \ar@<-0.3ex>[l]_{ {\scriptstyle #4}}}}

\title{A uniform open image theorem for \texorpdfstring{$\ell$-}-adic representations in positive characteristic}
\author{Emiliano Ambrosi}
\begin{abstract}
Let $k$ be a finitely generated field of characteristic $p > 0$ and $\ell$ a prime $\neq p$. Let $X$ be a smooth, separated, geometrically connected curve of finite type over $k$ and $\rho: \pi_1(X)\rightarrow GL_r(\Zl)$ a continuous representation of the \etale fundamental group of $X$ with image $\Pi$. Every $k$-rational point $x:Spec(k)\rightarrow X$ induces a local representation $\rho_x: \pi_1(Spec(k)) \rightarrow \pi_1(X) \rightarrow GL_r(\Zl)$ with image $\Pi_x$. The main result of this paper is that if every open subgroup of $\rho(\pi_1(X_{\overline k}))$ has finite abelianization, then the set $X_{\rho}^{ex}(k)$ of $k$-rational points such that $\Pi_x$ is not open in $\Pi$ is finite and there exists an integer $N\geq 1$, depending only on $\rho$, such that $[\Pi:\Pi_x]\leq N$ for all $x\in X(k)-X_{\rho}^{ex}(k)$.
This result can be applied to obtain uniform bounds for the $\ell$-primary torsion of group theoretic invariants in one dimensional families of varieties. For example, torsion of abelian varieties and Galois invariants of the geometric Brauer group. This extends to positive characteristic a previous result of Cadoret-Tamagawa in characteristic 0.
\end{abstract}
\maketitle
\tableofcontents
\section{Introduction}
\subsection{Notation}
In this paper $k$ is a field of characteristic $p>0$ with algebraic closure $k\subseteq \overline k$. A $k$-variety is reduced scheme separated and of finite type over $k$. For a $k$-variety $X$, write $|X|$ for the set of closed points and for every integer $d\geq 1$, $X(\leq d)$ for the set of all $x\in |X|$ with residue field $k(x)$ of degree $\leq d$ over $k$. If $d=1$ we often write $X(\leq 1)=X(k)$. Let $\ell$ be a prime always $\neq p$.
\subsection{Exceptional Locus}\label{inex}
From now on, let $X$ be a smooth geometrically connected $k$-variety. 
Let $\rho:\pi_1(X)\rightarrow GL_r(\Zl)$ be a continuous representation of the \etale fundamental group\footnote{As the choice of fibre functors will play no part in the following we will omit them for the notation for \etale fundamental group.} of $X$. By functoriality of the \etale fundamental group, every $x\in |X|$ induces a continuous group homomorphism $\pi_1(x)\rightarrow \pi_1(X)$, hence a “local" Galois\footnote{Recall that $\pi_1(x)\simeq \pi_1(Spec(k(x))$ identifies with the absolute Galois group of $k(x)$.} representation $\rho_{x}:\pi_1(x)\rightarrow \pi_1(X)\rightarrow GL_r(\Zl)$. Set
$$\Pi=\rho(\pi_1(X)) \quad  \Pi_{\overline k}=\rho(\pi_1(X_{\overline k}))\quad \Pi_x=\rho_{x}(\pi_1(x)).$$
Write $X_{\rho}^{gen}$ for the set of all $x\in |X|$ such that $\Pi_x\subseteq \Pi$ is an open subgroup of $\Pi$ and set
$$X_{\rho}^{ex}:=|X|-X_{\rho}^{gen}; \quad X_{\rho}^{gen}(\leq d):=X_{\rho}^{gen}\cap X(\leq d); \quad X_{\rho}^{ex}(\leq d):=X_{\rho}^{ex}\cap X(\leq d).$$ 
We call $X_{\rho}^{ex}$ the exceptional locus of $\rho$. The study of $X_{\rho}^{gen}(\leq d)$ is an important problem especially when the representation comes from a smooth proper morphism $f:Y\rightarrow X$ (see Subsection \ref{secapplication}), so that $\Pi_x$ controls fine arithmetic and geometric invariants of the family $Y_x$, $x\in |X|$. 
Since the Frattini subgroup of $\Pi$ is open (\cite[Pag. 148]{Serreweil}), a classical argument (\cite[Section 10.6]{Serreweil}) shows that if $k$ is Hilbertian (in particular if $k$ is finitely generated) there exists a $d\geq 1$ such that $X_{\rho}^{gen}(\leq d)$ is infinite.
\subsection{Uniform open image theorem}\label{statementUOI}
When $X$ is a curve and $k$ is finitely generated, one can go further, under a mild assumption on $\rho:\pi_1(X)\rightarrow GL_r(\Zl)$. 
\begin{definition}
A topological group $\Pi$ is Lie perfect\footnote{The terminology comes from the fact that if $\Pi$ is an $\ell$ adic Lie group this condition is equivalent to $Lie(\Pi)^{ab}=0$.} (or $LP$ for short) if every open subgroup of $\Pi$ has finite abelianization.
We say that $\rho:\pi_1(X)\rightarrow GL_r(\Zl)$ is Lie perfect (or $LP$ for short) if $\Pi$ is $LP$ and that $\rho$ is geometrically Lie perfect (or $GLP$ for short) if $\Pi_{\overline k}$ is $LP$.
\end{definition}
With this terminology, we can state our main result, which is an extension of \cite[Theorem 1.1]{UOI1} to positive characteristic.
\begin{theorem}\label{UOI}
Assume that $X$ is a curve and $k$ is finitely generated. If $\rho$ is $GLP$, then $X_{\rho}^{ex}(\leq 1)$ is finite and there exists an integer $N\geq 1$, depending only on $\rho$, such that $[\Pi:\Pi_x]\leq N$ for all $x\in X_{\rho}^{gen}(\leq 1)$.
\end{theorem}
When $X(k)$ is infinite, Theorem \ref{UOI} gives us uniform boundedness results that are impossible to achieve using the arguments in \cite[Section 10.6]{Serreweil}; see for example Corollaries \ref{uniformintroab} and \ref{uniformintrok3}.
\subsection{Strategy}
While the general strategy of the proof of Theorem \ref{genus} is similar to the one of \cite[Theorem 1.1]{UOI1}, the technical details are more complicated in positive characteristic. Indeed, the proof of Theorem \ref{UOI} is based on the genus computations, via the Riemann-Hurwitz formula, of careful chosen abstract modular curves. In positive characteristic, the Riemann-Hurwitz formula involves wild inertia terms and - even assuming $\ell\neq p$ - controlling those wild inertia terms is rather delicate. To deal with them, we generalize the computations made in \cite{UB}.
\subsubsection{Abstract modular scheme}\label{ingenusgonality}
For every open subgroup $U\subseteq \Pi$ write $f_U:X_U\rightarrow X$ for the connected \etale cover corresponding to the open subgroup $\rho^{-1}(U)\subseteq \pi_1(X)$ and $k_U$ for the smallest separable field extension of $k$ over which $X_U$ is geometrically connected. Write $U_{\overline k}=U\cap \Pi_{\overline k}$ and recall the following anabelian dictionary.
\begin{fact}\label{formalprop}
For every open subgroup $U\subseteq \Pi$ the following hold:
\begin{enumerate}
\item For every $x\in |X|$, we have that $\Pi_x\subseteq U$ if and only if $x$ lifts to a $k(x)$-rational point on $X_U$;
\item The cover $X_{U_{\overline k}}\rightarrow X_{\overline k}$ corresponding to the open subgroup $U_{\overline k}\subseteq \Pi_{\overline k}$ is $X_U\times_{k_U}\overline k\rightarrow X_{\overline k}$.
\end{enumerate}
\end{fact}
In view of Fact \ref{formalprop}, we call $X_U$ the connected abstract modular scheme associated to $U$. Fact \ref{formalprop} enabled Cadoret-Tamagawa in \cite{UOI1} to construct a projective system of
abstract modular schemes (whose definition is recalled in Section \ref{definitionprojective}):
$$f_{n}:\mathcal X_{n}:=\coprod_{U\in \mathcal C_{n}(\Pi)}X_U\rightarrow X.$$
This system has the property that if $x\in |X|$ does not lift to a $k(x)$-rational point of $\mathcal X_{n}$ for some $n\geq 1$, then $\Pi_x\subseteq\Pi$ is not an open subgroup; see Lemma \ref{lemmafondamentale}. 
The finiteness of $X^{ex}_{\rho}(\leq d)$ can be then formulated in diophantine terms as follows:
\begin{equation*}\label{equationintrodef}
\text{(1): The image of $\varprojlim \mathcal X_n(\leq d)\rightarrow X(\leq d)$ is finite.}
\end{equation*}
To prove (1) it is enough to show
\begin{equation*}\label{equationintrodef2}
\text{(2): The set $\mathcal X_n(\leq d)$ is finite for $n\gg 0$.}
\end{equation*}
\subsubsection{Growth of Genus}
If $X$ is a curve, $k$ is finitely generated and $d=1$, by \cite{samuel} and an argument of Voloch (see \cite[Theorem 3]{mordel} for more details), the finiteness of $\mathcal X_{n}(k)$ is controlled by the genus $g_U$ of the smooth compactification of $X_{U_{\overline k}}$ for $U\in\mathcal C_n(\Pi)$.
\begin{fact}\label{mordel}
If $k$ is finitely generated of positive characteristic, there exists an integer $g\geq 2$, depending only on $k$, such that for every smooth proper $k$-curve $Y$ with genus $\geq g$, the set $Y(k)$ is finite.
\end{fact}
Fact \ref{mordel} reduces $(2)$ to the geometric Theorem \ref{genus} below, which extends  \cite[Theorem 3.4]{UOI1} to positive characteristic. Write $\Pi_{\overline k}(n):=\Ker(\Pi_{\overline k}\rightarrow GL_r( \Zl/\ell^n))$.
\begin{theorem}\label{genus}
Assume that $X$ is a curve, $\rho$ is $GLP$ and $\ell\neq p$. Then for every closed but not open subgroup $C\subseteq \Pi_{\overline k}$ we have $$\lim_{n\to +\infty}g_{C\Pi_{\overline k}(n)}=+\infty.$$
\end{theorem}
To prove Theorem \ref{genus}, one may assume $k=\overline k$, hence that $\Pi=\Pi_{\overline k}$. We first replace $X_{C\Pi(n)}\rightarrow X$ with a Galois cover $X_{\widetilde \Pi_C(n)}\rightarrow X$, closely related to the Galois closure of $X_{C\Pi(n)}\rightarrow X$, and we use the $GLP$ hypothesis to show that the genus of $X_{\widetilde \Pi_C(n)}$ goes to infinity. Then we translate into group theoretical terms the Riemann-Hurwitz formula for $X_{\widetilde \Pi_C(n)} \rightarrow  X_{C\Pi(n)}$ to
show that the genus of $X_{\widetilde \Pi_C(n)}$ tends to infinity (if and) only if the genus of $X_{C\Pi(n)}$ does. Here, we use crucially that $\ell\neq p$ to control the wild inertia terms appearing in the Riemann-Hurwitz formula for $X_{\widetilde \Pi_C(n)} \rightarrow  X_{C\Pi(n)}$. This part of the argument is significantly more difficult than in the proof of \cite[Theorem 3.4]{UOI1}.
\subsection{Applications to motivic representations}\label{secapplication}
Let $f:Y\rightarrow X$ be a smooth proper morphism and let $\ell\neq p$ be a prime. For $x\in X$, choose a geometric point $\overline x$ over $x$ and set $Y_x$ (resp. $Y_{\overline x}$) for the fibre of $f$ at $x$ (resp. $\overline x$). By smooth proper base change $R^if_*\Zl(j)$ is a lisse sheaf hence, for every $x\in |X|$, gives rise to a continuous representation 
$$\rho_{\ell}:\pi_1(X)\rightarrow GL(H^i(Y_{\overline x},\Zl(j)))$$
such that $\rho_{\ell, x}:\pi_1(x)\rightarrow GL(H^i(Y_{\overline x},\Zl(j)))$ identifies with the natural Galois action of $\pi_1(x)$ on $H^i(Y_{\overline x},\Zl(j))$. By \cite[Theorem 5.8]{UOI1}, the representation $\rho_{\ell}$ is $GLP$, so that we can apply Theorem \ref{UOI} to it.
\subsubsection{Uniform boundedness \texorpdfstring{$\ell$-}-primary torsion of abelian schemes}
Let $f:Y\rightarrow X$ be a $g$-dimensional abelian scheme. For $x\in X$ and any integer $n\geq 1$, write $Y_{\overline{x}}[\ell^{n}]:Y_{\overline{x}}[\ell^{n}](\overline{k(x)})$ for the $\ell^n$-torsion of $Y_{\overline x}$ and set
$$Y_{\overline x}[\ell^{\infty}]:=\bigcup_n Y_{\overline x}[\ell^{n}];\quad T_{\ell}(Y_{\overline x}):=\varprojlim_nY_{\overline x}[\ell^{n}].$$
Since $k$ is finitely generated, $Y_x[\ell^{\infty}](k(x))(=Y_{\overline x}[\ell^{\infty}]^{\pi_1(x)})$ is finite by the Mordell-Weil theorem. From the $\pi_1(x)$-equivariant isomorphisms
$$T_{\ell}(Y_{\overline x})\simeq H^{2g-1}(Y_{\overline x},\Zl(g));\quad T_{\ell}(Y_{\overline x})\otimes \Ql/\Zl\simeq Y_x[\ell^{\infty}](\overline{k(x)})$$
and Theorem \ref{UOI}, we obtain the following uniform bound for $Y_{x}[\ell^{\infty}](k)$, $x \in X(k)$.
\begin{corollary}\label{uniformintroab}
Assume that $k$ is finitely generated, $X$ is a curve and $f:Y\rightarrow X$ is an abelian scheme. There exists an integer $N\geq 1$, depending only on $f:Y\rightarrow X$ and $\ell$, such that $|Y_x[\ell^{\infty}](k)|\leq N$ for every $x\in X(k)$.
\end{corollary}
\proof
Since $\Pi_{\ell}:=\rho_{\ell}(\pi_1(X))$ is a compact $\ell$-adic Lie group, it is topologically finitely generated hence it has only finitely many open subgroups of bounded index. So, by Theorem \ref{UOI}, the set of subgroups $\Pi_{\ell,x}\subseteq \Pi_{\ell}$ appearing as $\rho_{\ell,x}(\pi_1(x))$ for $x\in X(k)$ is finite. In particular, the set of abelian groups $\{Y_{\overline x}[\ell^{\infty}]^{\pi_1(x)}\simeq Y_x[\ell^{\infty}](k) \text{ } | \text{ } x\in X(k)\}$ is finite.\endproof
\subsubsection{Further applications}
In the subsequent paper \cite{mioneron}, Theorem \ref{UOI} is used to prove the following results. For $x\in |X|$, let $Br(Y_{\overline x})^{\pi_1(x)}[\ell^{\infty}]$ denote the Galois invariants of the $\ell$-primary torsion of the geometric Brauer group $Br(Y_{\overline x}):=H^2(Y_{\overline x},\Gm)$ of $Y_x$.
\begin{corollary}\label{uniformintrok3}
Assume that $k$ is finitely generated and that $X$ is a curve with generic point $\eta$. Then
\begin{itemize}
\item \cite[Corollary 1.6.3.1]{mioneron}: Assume that all the closed fibres of $f:Y\rightarrow X$ satisfy\footnote{This holds, for example, if $f:Y\rightarrow X$ is a family of abelian varieties or of K3 surfaces.} the $\ell$-adic Tate conjecture for divisors (\cite{tateconj}). Then there exists an integer $N\geq 1$, depending only on $f:Y\rightarrow X$ and $\ell$, such that $|Br(Y_{\overline x})^{\pi_1(x)}[\ell^{\infty}]|\leq N$ for every $x\in X(k)$.
\item \cite[Corollary 1.6.1.2]{mioneron}: For all but at most finitely many $x\in X(k)$, the rank of the N\'eron-Severi group of $Y_{\overline x}$ is the same as the one of the N\'eron-Severi group of $Y_{\overline \eta}$
\end{itemize}
\end{corollary}
Corollaries \ref{uniformintroab} and \ref{uniformintrok3} are extensions to positive characteristic of previous results obtained in \cite{UOI1}, \cite[Thm. 1.6, Cor. 1.7]{VAV} and \cite{brauer}.
\subsection{Organization of the paper}
The paper is organized as follows. In Section \ref{secgenus} we prove Theorem \ref{genus}. In Section \ref{secmain} we recall the construction of a projective system of abstract modular schemes $\mathcal X_n\rightarrow X$, parametrizing points with small image and some facts about them. After this, we prove Theorem \ref{UOI}. In Subsection \ref{secquestion1}, we discuss possible extensions of Theorem \ref{UOI} to points of bounded degree.
All the results and the proofs in this paper work in the characteristic zero setting but, since this situation is already
treated in \cite{UOI1}, we will assume that $p>0$ to simplify the exposition.
\subsection{Acknowledgement} 
This paper is part of the author's Ph.D. thesis under the supervision of Anna Cadoret. He thanks her for suggesting the problem and the appropriate references. He is also very grateful for her careful re-readings of this paper and her constructive suggestions.
%\section{Notation and preliminaries}\label{secnotation}
%
%\subsection{Ramification}\label{notramification}
%\numberwithin{equation}{subsection} 
%We use the capital letters $K,K',L$ to denote local fields with algebraically closed residue fields.
%If $f:Z\rightarrow Y$ is a (possibly ramified) Galois cover of smooth proper curves over an algebraically closed field $k$, define the ramification group at $z\in |Z|$ as the ramification group of the Galois extension of local fields $k((f(z)))\rightarrow k((z))$. 
%If $y\in |Y|$ we call the ramification group $\Pi_y$ of $f$ over $y$ any of the ramification group at $z$, for any $z$ with $f(z)=y$. Since the cover is Galois the conjugacy class of the ramification group does not depend on the choice of $z$.
\section{Proof of Theorem \texorpdfstring{ \ref{genus} }  -}\label{secgenus}
\subsection{Notation}
\subsubsection{}
For a group $\Gamma$ and subgroups $I,H\subseteq \Gamma$ write 
$$K_H(\Gamma):=\cap_{g\in \Gamma}gHg^{-1} \quad \text{and} \quad I_H:=I/(I\cap K_H(\Gamma))$$ for the largest normal subgroup of $\Gamma$ contained in $H$ and the largest quotient of $I$ that acts faithfully on $\Gamma/H$.
For every closed subgroup $\Gamma \subseteq GL_r(\Zl)$, write $\Gamma(n):=\Ker(\Gamma\rightarrow GL_r(\Zl/\ell^n))$ and $\Gamma_n:=\Image(\Gamma\rightarrow GL_r(\Zl/\ell^n)).$
We use $\twoheadrightarrow$ and $\hookrightarrow$ to denote surjective and injective maps respectively.
\subsubsection{}
From now on we retain the notation of Theorem \ref{genus}. By Fact \ref{formalprop}(2), we may assume $k=\overline k$, hence that $\Pi=\Pi_{\overline k}$ is LP.
Set 
$$\widetilde \Pi_C(n):=\Ker(\Pi_C\twoheadrightarrow (\Pi_n)_{C_n}) \quad \text{and} \quad \Pi_C(n):=\Ker(\Pi\twoheadrightarrow (\Pi_n)_{C_n}).$$
The following exact diagram summarizes the situation:
\begin{center}
\begin{tikzcd}[column sep=tiny, row sep=tiny]
 &  &  & 1\arrow{d}\\

& 1\arrow{rd} & & \widetilde \Pi_C(n)\arrow{dd}\arrow[bend right=50, two heads]{ddll}\arrow[bend left=50, two heads]{ddrr} & & 1\arrow{ld}\\

1\arrow{rd} & & K_C(\Pi)\arrow[hook]{ur}\arrow{dr}\arrow{ld}& & \Pi(n)\arrow[hook]{ul}\arrow{ld}\arrow[two heads]{dr} & & 1\arrow{ld}\\

& K_{C_n}(\Pi_n)\arrow{rd} & & \Pi\arrow{dd}\arrow{dl}\arrow{rd} & & \Pi_C(n)\arrow{dl}\\

& & \Pi_n\arrow{dr}\arrow{dl} & & \Pi_C\arrow{dl}\arrow{dr}\\

& 1 & & (\Pi_n)_{C_n}\arrow{d}\arrow{ld}\arrow{rd}& & 1\\

& & 1 & 1 & 1.

\end{tikzcd}
\end{center}
Since if $g\in \Pi$ acts trivially on $\Pi_n/C_n$ for every $n$ then it acts trivially on $\Pi/C$, the collection $\{\Pi_C(n)\}_{n\in \N}$ is a fundamental system of neighbourhoods of the identity in $\Pi_C$. \subsection{Preliminary reduction}\label{preliminaryyyyy}
In this section we show that we can assume that for every integer $n\geq 1$:
\begin{enumerate}
\item $K_C(\Pi)=K_C(C\Pi(n))$;
\item $\widetilde \Pi_C(1)/\widetilde \Pi_C(n)$ is an $\ell$-group.
\end{enumerate}
Since we are interested in the asymptotic behaviour of $g_{C\Pi(n)}$ we can freely replace $\Pi$ with $C\Pi(n_0)$ for some integer $n_0\geq 1$. So:
\begin{enumerate}
\item Follows from the fact the increasing sequence $K_C(\Pi)\subseteq K_C(C\Pi(1))\subseteq ...\subseteq K_C(C\Pi(n))\subseteq ...$ of closed subgroups of $\Pi$ stabilizes (\cite[Theorem 6.1]{UOI1});
\item Follows if we prove that $\widetilde \Pi_C(n_0)/\widetilde \Pi_C(n)$ is an $\ell$-group for some integer $n_0\geq 1$ and any $n>n_0$. Write $A_n:=\widetilde \Pi_C(n)/\Pi(n)$. Using the commutative exact diagram
\begin{center}
\begin{tikzcd}[column sep=scriptsize, row sep=scriptsize]
1 \arrow{r} &  \Pi(n)\arrow{r}\arrow[hook]{d} & \widetilde \Pi_C(n)\arrow{r}\arrow[hook]{d} & A_n\arrow{r}\arrow{d} & 1 \\
1 \arrow{r} & \Pi(1)\arrow{r} & \widetilde \Pi_C(1)\arrow{r} & A_1\arrow{r} & 1 
\end{tikzcd}
\end{center}
we find an exact sequence
$$1\rightarrow B_{\ell,n}\rightarrow \widetilde \Pi_C(1)/\widetilde \Pi_C(n)\rightarrow A_1/A_n\rightarrow 1,$$
where $B_{\ell,n}$ is a quotient of $\Pi(1)/\Pi(n)$, hence an $\ell$-group.
Since $A_1$ is finite, for some $n_0\gg 0$ and any $n\geq n_0$ the surjection $A_1/A_n\twoheadrightarrow A_1/A_{n-1}$ is an isomorphism. 
The (non abelian) snake lemma applied to the commutative diagram
\begin{center}
\begin{tikzcd}[column sep=scriptsize, row sep=scriptsize]
1 \arrow{r} & B_{\ell,n}\arrow{r}\arrow{d} & \widetilde \Pi_C(1)/\widetilde \Pi_C(n)\arrow{r}\arrow[two heads]{d} & A_1/A_n\arrow{r}\arrow{d}{\simeq} & 1 \\
1 \arrow{r} & B_{\ell,n-1}\arrow{r} & \widetilde \Pi_C(1)/\widetilde \Pi_C(n-1)\arrow{r} & A_1/A_{n-1}\arrow{r} & 1,
\end{tikzcd}
\end{center}
shows that 
$$\widetilde \Pi_C(n-1)/\widetilde \Pi_C(n)=\Ker(\widetilde \Pi_C(1)/\widetilde \Pi_C(n)\rightarrow \widetilde \Pi_C(1)/\widetilde \Pi_C(n-1))\subseteq  B_{\ell,n},$$  
hence that $\widetilde \Pi_C(n-1)/\widetilde \Pi_C(n)$ is an $\ell$-group.
We conclude by induction on $n\geq n_0$ using the exact sequence 
$$1\rightarrow \widetilde \Pi_C(n-1)/\widetilde \Pi_C(n)\rightarrow \widetilde \Pi_C(n_0)/\widetilde \Pi_C(n)\rightarrow \widetilde \Pi_C(n_0)/\widetilde \Pi_C(n-1)\rightarrow 1.$$
\end{enumerate} 
So, from now on we may and do assume that (1) and (2) hold.
\subsection{\texorpdfstring{ $\mathbf{g}_{\widetilde \Pi_C(n)}\to +\infty}-$}\label{usingGLP}
We use that $\Pi$ is Lie perfect and $X_{\widetilde \Pi_C(n)}\rightarrow X$ is Galois. Since $C$ is closed but not open in $\Pi$, $|(\Pi_n)_{C_n}|\to +\infty$ hence $g_{\widetilde \Pi_C(n)}\to \infty$ as soon as $\sup{g_{\widetilde \Pi_C(n)}}>1$. Indeed, assume that $g_{\widetilde \Pi_C(n_0)}>1$ for some $n_0\geq 1$. Then, for every $n>n_0$, the Riemann Hurwitz formula for $X_{\widetilde \Pi_C(n)}\rightarrow X_{\widetilde \Pi_C(n_0)}$  yields
$$\lim_{n\to +\infty}2g_{\widetilde \Pi_C(n)}-2\geq\lim_{n\to +\infty}\frac{|(\Pi_{n})_{C_{n}}|}{|(\Pi_{n_0})_{C_{n_0}}|}(2g_{\widetilde \Pi_C(n_0)}-2)= +\infty.$$ 
So it remains to show that $\sup{g_{\widetilde \Pi_C(n)}}=1$ and  $\sup{g_{\widetilde \Pi_C(n)}}=0$ are not possible.
\subsubsection{\texorpdfstring{ $\sup{g_{\widetilde \Pi_C(n)}}=1$}-}\label{sectiongenus1}
Assume $\sup{g_{\widetilde \Pi_C(n)}}=1$. Then there exists $n_0$ such that for all $n\geq n_0$ the smooth compactification of $X_{\widetilde \Pi_C(n)}$ is an elliptic curve. Since finite morphisms between elliptic curves are unramified, the Galois group $\widetilde \Pi (n_0)/\widetilde \Pi (n)\simeq\Pi_C(n_0)/\Pi_C(n)$ of $X_{\widetilde \Pi_C(n)}\rightarrow X_{\widetilde \Pi_C(n_0)}$ would be a quotient of the \etale fundamental group of the smooth compactification of $X_{\widetilde \Pi_C(n_0)}$. In particular it would be abelian, hence $$\Pi_C(n_0)=\varprojlim_n \Pi_C(n_0)/\Pi_C(n)$$ would be abelian and infinite. But this contradicts the fact that $\Pi$ is Lie perfect, since $\Pi_C(n_0)$ would be an infinite abelian quotient of the open subgroup $\widetilde \Pi_C(n_0)$ of $\Pi$.
\subsubsection{\texorpdfstring{ $\sup{g_{\widetilde \Pi_C(n)}}=0$}-}
Assume $\sup{g_{\widetilde \Pi_C(n)}}=0$. This means that for all $n\geq 0$, the smooth compactification of $X_{\widetilde \Pi_C(n)}$ is isomorphic to $\PP^1$. So the Galois group $\widetilde \Pi_C(1)/\widetilde \Pi_C(n)\simeq \Pi_C(1)/\Pi_C(n)$ of $X_{\widetilde \Pi_C(n)}\rightarrow X_{\widetilde \Pi_C(1)}$ is a subgroup of $PGL_2(k)$. We use the following:
\begin{fact}\cite[Corollary 10]{notes}\label{classification}
Suppose that $k$ is an algebraically closed field of characteristic $p> 0$. A finite subgroup of $PGL_2(k)$ is isomorphic to one of the following groups:
\begin{itemize}
\item A cyclic group;
\item A dihedral group $D_{2m}$ of order $2m$, for some $m>0$;
\item $A_4,A_5,S_4$;
\item An extension $1\rightarrow A\rightarrow \Pi\rightarrow Q\rightarrow 1$, with $A$ an elementary abelian p-group and $Q$ a cyclic group of prime-to-$p$ order;
\item $PSL_2(\F_{p^r})$, for some $r>0$;
\item $PGL_2(\F_{p^r})$, for some $r>0$;
\end{itemize}
where $\F_{p^r}$ denotes the finite field with $p^r$ elements.
\end{fact}
From Fact \ref{classification} and the fact that $\Pi_C(1)/\Pi_C(n)$ is an $\ell$-groups by Section \ref{preliminaryyyyy}(2), the only possibility is that $\Pi_C(1)/\Pi_C(n)$ is a cyclic group or $\ell=2$ and $\Pi_C(1)/\Pi_C(n)\simeq D_{2^m}$. If the groups 
$\Pi_C(1)/\Pi_C(n)$ are abelian for $n\gg 0$ we can conclude as in \ref{sectiongenus1}. So assume $\ell=2$ and $\Pi_C(1)/\Pi_C(n)\simeq D_{2^m}$. Since $D_{2^m}$  fits into an exact sequence
$$0\rightarrow \Z/2^{m-1} \rightarrow D_{2^m}\rightarrow \Z/2\times \Z/2\rightarrow 0,$$
the exactness of $\varprojlim$ on finite groups yields an infinite abelian open subgroup $\Z_2\subseteq\displaystyle  \varprojlim_n \Pi_C(1)/\Pi_C(n)\simeq \Pi_C(1)$, 
and we conclude as in \ref{sectiongenus1}.
\subsection{Definition of \texorpdfstring{$\lambda$}-}\label{subsectiondefinitionoflamba}
If $f:Y\rightarrow X$ is a cover we define $$\lambda_{Y/X}:=\frac{2g_Y-2}{deg(f)}.$$ 
The following directly follows from the Riemann-Hurwitz formula.
\begin{lemma}\label{lemmaRH}
Let $... \rightarrow X_{n+1}\rightarrow X_n\rightarrow ...\rightarrow ...\rightarrow X$ be a sequence of finite covers of smooth proper connected curves over an algebraically closed field $k$. Then $\lambda_{X_{n+1}/X}\geq \lambda_{X_{n}/X}$. Assume furthermore that $Deg(X_n\rightarrow X)\to +\infty$. Then $g_{X_n}\to +\infty$ if and only if $\displaystyle \lim_{n\to +\infty}\lambda_{X_n/X}>0$

%\underset{n\to +\infty}{\lim} \lambda_{X_n/X}>0$.
\end{lemma}
For an open subgroup $U\subseteq \Pi$ write $\lambda_U:=\lambda_{X_U/X}$. With this notation, applying Lemma \ref{lemmaRH} to
\begin{center}
\begin{tikzcd}[column sep=scriptsize, row sep=scriptsize]
...\arrow{r} &X_{\widetilde \Pi_C(n+1)}\arrow{r}\arrow{d}&X_{\widetilde \Pi_C(n)}\arrow{r}\arrow{d}&... \arrow{r} & X_{\widetilde \Pi_C(1)}\arrow{r}\arrow{d}& X\arrow[equal]{d}\\
...\arrow{r} &X_{C\Pi(n+1)}\arrow{r}&X_{C\Pi(n)}\arrow{r}&... \arrow{r} & X_{C\Pi(1)}\arrow{r}& X\\
\end{tikzcd}
\end{center}
one gets inequalities
\begin{center}
\begin{tikzcd}[column sep=tiny, row sep=tiny]
\lambda_{C\Pi(n+1)}\arrow[phantom]{r}{\geq} & \lambda_{C\Pi(n)}\\
\lambda_{\widetilde \Pi_C(n+1)}\arrow[phantom]{r}{\geq}\arrow[phantom]{u}{\rotatebox{90}{$\geq $}} &\lambda_{\widetilde \Pi_C(n)}\arrow[phantom]{u}{\rotatebox{90}{$\geq $}}
\end{tikzcd}
\end{center}
hence $\lambda_{\widetilde \Pi_C}:=\underset{n\to +\infty}{\lim} \widetilde \Pi_C(n)$ and $\lambda_{C}:=\underset{n\to +\infty}{\lim} \lambda_{C\Pi(n)}$ exist with $\lambda_{\widetilde \Pi_C}\geq \lambda_{C}$.
Also, since $C\subseteq \Pi$ is closed but not open
\begin{enumerate}
\item $|(\Pi_n)_{C_n}|\to +\infty$, hence $g_{\widetilde \Pi_{C}(n)}\to +\infty$ if and only if $\lambda_{\widetilde \Pi_C}>0$;
\item $|\Pi_n/C_n|\to +\infty$, hence $g_{C\Pi(n)}\to +\infty$ if and only if $\lambda_{C}>0$.
\end{enumerate}
By Section \ref{usingGLP}, $\lambda_{\widetilde \Pi_C}>0$ hence it is enough to show that $\lambda_{\widetilde \Pi_C}=\lambda_C$. 
The remaining part of this section is devoted to the proof of this fact. 
\subsection{Inertia subgroups}\label{subsectioncomparisonUOI}
Consider the commutative diagram:
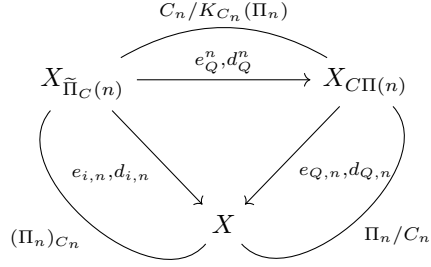
\begin{center}
\begin{tikzcd}
X_{\widetilde \Pi_C(n)}\arrow[bend left, dash]{rr}{C_n/K_{C_n}(\Pi_n)}\arrow{rr}{e^n_Q,d^n_Q}\arrow{ddr}[swap]{e_{i,n}, d_{i,n}}
\arrow[bend right=90, dash, swap]{ddr}{(\Pi_n)_{C_n}} & & X_{C\Pi(n)} \arrow{ddl}{e_{Q,n}, d_{Q,n}}\arrow[bend left=90, dash]{ddl}{\Pi_n/C_n}\\\\
& X
\end{tikzcd}
\end{center}
Suppose that $X^{cpt}-X=\{P_1,...,P_r\}$ and denote with $I_i\subseteq \Pi$ the image via $\pi_1(X)\twoheadrightarrow \Pi$ of the inertia group of the point $P_i$. The situation is then the following.
\begin{itemize}
\item $X_{\widetilde \Pi_C(n)}\rightarrow X$ is a Galois cover with Galois group $(\Pi_{n})_{C_{n}}$. The inertia group and the ramification index of any point of $X_{\widetilde \Pi_C(n)}$ over $P_i$ are given\footnote{Since the cover is Galois the conjugacy class of the ramification group does not depend on the choice of the point over $P_i$} by $(I_{i,n})_{C_n}\subseteq (\Pi_{n})_{C_{n}}$ and $e_{i,n}:=|(I_{i,n})_{C_n}|$. Write $((I_{i,n})_{C_n})_j\subseteq (\Pi_{n})_{C_{n}}$ for the $j^{th}$-ramification group in lower numbering (see \cite[Section 1, IV]{Serre}) over the point $P_i$ and $(e_{i,n})_j$ for its cardinality. Finally set $d_{i,n}$ for the exponent of the different of any point of $X_{\widetilde \Pi_C(n)}$ over $P_i$.
\item $X_{C \Pi(n)}\rightarrow X$ is the cover corresponding to the open subgroup $C\Pi(n)\subseteq \Pi$. If $Q\in X_{C\Pi(n)}$ is over $P_i$ we denote with $e_{Q,n}, d_{Q,n}$ the ramification index and the exponent of the different of $Q$ over $P_i$.
\item $X_{\widetilde \Pi_C(n)}\rightarrow X_{C\Pi(n)}$ is a Galois cover with Galois group $C_n/K_{C_n}(\Pi_n)\subseteq (\Pi_{n})_{C_{n}}$ and there is a natural bijection of sets $$\{Q\in X_{C\Pi(n)} \text{ }| \text{ } Q|P_i\}\simeq (I_{i,n})_{C_n}\backslash \Pi_n/C_n.$$  
If $Q$ correspond to the element $(I_{i,n})_{C_n}x\in  (I_{i,n})_{C_n}\backslash \Pi_n/C_n$, then the inertia group and the ramification index at $Q$ are given by $Stab_{(I_{i,n})_{C_{n}}}((I_{i,n})_{C_n}x)$ and $|Stab_{(I_{i,n})_{C_{n}}}((I_{i,n})_{C_n}x)|:=e^n_Q$. The $j^{th}$-ramification group is given by $((I_{i,n})_{C_n})_j\cap Stab_{(I_{i,n})_{C_{n}}}(x)$. Write $|((I_{i,n})_{C_n})_j\cap Stab_{(I_{i,n})_{C_{n}}}(x)|=(e^n_Q)_j$.
\end{itemize}
By \cite[Section 4, III, Pag. 51]{Serre} we have the following relations:
$$e_{i,n}=e^n_Qe_{Q,n};\quad d_{i,n}=d^n_Q+e^n_Qd_{Q,n}; \quad \sum_{Q|P_i}e_{Q,n}=|\Pi_n/C_n|.$$
\subsection{Comparison}
Using the Riemann-Hurwitz formula we get
$$\lambda_{C\Pi(n)}=2g_X-2+\frac{1}{|\Pi_n/C_n|}\sum_{1\leq i\leq r}\sum_{Q|P_i}d_{Q,n} \quad \text{and} \quad \lambda_{\widetilde \Pi_C(n)}=2g_X-2+\sum_{1\leq i\leq r}\frac{d_{i,n}}{e_{i,n}},$$
hence $$\lambda_{\widetilde \Pi_C(n)}-\lambda_{C\Pi(n)}=\frac{1}{|\Pi_n/C_n|}\bigg (\sum_{1\leq i\leq r}\frac{d_{i,n}|\Pi_n/C_n|}{e_{i,n}}-\sum_{Q|P_i}\frac{d_{Q,n}e_{i,n}}{e_{i,n}}\bigg )$$
$$=\frac{1}{|\Pi_n/C_n|}\bigg (\sum_{1\leq i\leq r}\frac{d_{i,n}\sum_{Q|P_i}e_{Q,n}}{e_{i,n}}-\sum_{Q|P_i}\frac{d_{Q,n}e_{i,n}}{e_{i,n}}\bigg )=\frac{1}{|\Pi_n/C_n|}\bigg (\sum_{1\leq i\leq r}\sum_{Q|P_i}\frac{d_{i,n}e_{Q,n}-d_{Q,n}e_{i,n}}{e_{i,n}}\bigg )$$
$$=\frac{1}{|\Pi_n/C_n|}\bigg (\sum_{1\leq i\leq r}\sum_{Q|P_i}\frac{d_{i,n}-d_{n,Q}e^n_Q}{e^n_Q}\bigg )=\frac{1}{|\Pi_n/C_n|}\bigg (\sum_{1\leq i\leq r}\sum_{Q|P_i}\frac{d^n_Q}{e^n_Q}\bigg ).$$
So it is enough to show that for every integer $1\leq i\leq r$ one has
 $$\lim_{n\to \infty}\frac{1}{|\Pi_n/C_n|}\sum_{Q|P_i}\frac{d^n_Q}{e^n_Q}=0.$$
By \cite[Proposition 4, IV, Pag. 64]{Serre} we have $$\frac{1}{|\Pi_n/C_n|}\sum_{Q|P_i}\frac{d^n_Q}{e^n_Q}=\frac{1}{|\Pi_n/C_n|}\sum_{j\geq 0}\sum_{Q|P_i}\frac{(e^n_Q)_j-1}{e^n_Q}.$$
\subsection{Galois formalism}\label{Galoisformalism}
Consider the surjection $$\phi_j:((I_{i,n})_{C_n})_j\backslash \Pi_n/C_n\rightarrow (I_{i,n})_{C_n}\backslash \Pi_n/C_n$$ and recall the following elementary lemma.
\begin{lemma}[{\cite[Lemma 4.3]{UB}}]\label{lemmaditeoriadeigruppistupido}
Let $G$ be a finite group and $H\subseteq G$ a normal subgroup. Let $X$ be a finite set on which $G$ acts and consider the natural surjection $q:H\backslash X\twoheadrightarrow G\backslash X$. If $Gx\in G\backslash X$ then
$$|q^{-1}(Gx)|=\frac{|G||Stab_{G}(Gx)\cap H|}{|H||Stab_{G}(Gx)|}.$$
\end{lemma}
If, under the bijection $$(I_{i,n})_{C_n}\backslash \Pi_n/C_n\simeq \{Q\in X_{C\Pi(n)} \text{ }| \text{ } Q|P_i\},$$ 
the element $(I_{i,n})_{C_n}x\in  (I_{i,n})_{C_n}\backslash \Pi_n/C_n$ corresponds to the point $Q\in X_{C\Pi(n)}$ above $P_i$, by Lemma \ref{lemmaditeoriadeigruppistupido} we have
$$|\phi_j^{-1}((I_{i,n})_{C_n}x)|=\frac{|(I_{i,n})_{C_n}|}{|((I_{i,n})_{C_n})_j|}\frac{|((I_{i,n})_{C_n})_j\cap Stab_{(I_{i,n})_{C_{n}}}((I_{i,n})_{C_n}x)|}{|Stab_{(I_{i,n})_{C_{n}}}((I_{i,n})_{C_n}x)|}=\frac{e_{i,n}}{(e_{i,n})_j}\frac{(e^n_Q)_j}{e^n_Q}.$$
Summing over all the $Q\in X_{C\Pi(n)}$  above $P_i$, we get $$\sum_{Q|P_i}\frac{e_{i,n}}{(e_{i,n})_j}\frac{(e^n_Q)_j}{(e^n_Q)}=|((I_{i,n})_{C_n})_j\backslash \Pi_n/C_n|.$$
A similar reasoning gives $$\sum_{Q|P_i}\frac{e_{i,n}}{e^n_Q}=|\Pi_n/C_n|,$$ hence 
$$\frac{1}{|\Pi_n/C_n|}\sum_{j\geq 0}\sum_{Q|P_i}\frac{(e^n_Q)_j-1}{e^n_Q}=\frac{1}{|\Pi_n/C_n|}\bigg (\sum_{Q|P_i}\frac{e^n_Q-1}{e^n_Q}\frac{e_{i,n}}{e_{i,n}}\bigg )+\frac{1}{|\Pi_n/C_n|}\bigg (\sum_{j\geq 1}\sum_{Q|P_i}\frac{(e^n_Q)_j-1}{e^n_Q}\frac{(e_{i,n})_je_{i,n}}{(e_{i,n})_je_{i,n}}\bigg ).$$
The first term is 
$$\frac{1}{|\Pi_n/C_n|}\sum_{Q|P_i}1-\frac{1}{|\Pi_n/C_n|}\frac{1}{e_{i,n}}\sum_{Q|P_i}\frac{e_{i,n}}{e^n_Q}=\frac{|(I_{i,n})_{C_n}\backslash \Pi_n/C_n|}{|\Pi_n/C_n|}-\frac{1}{|(I_{i,n})_{C_n}|}.$$ 
Recall the following:
\begin{fact}[{\cite[Theorem 2.1]{UOI1}}]\label{limite}
Let $\Pi\subseteq GL_r(\Zl)$ be a closed subgroup and $C\subseteq \Pi$ a closed but not open subgroup. If $K_C(\Pi)=K_C(C\Pi(n))$ for every integer $n\geq 0$, then for every closed subgroup $I\subseteq \Pi$ one has
$$\lim_{n\to +\infty} \frac{|I_n\backslash \Pi_n/C_n|}{|\Pi_n/C_n|}=\frac{1}{|I_C|}.$$
\end{fact}
Since $\varprojlim(I_{i,n})_{C_n}=(I_i)_C,$ Fact \ref{limite} and Section \ref{preliminaryyyyy}(1) show that 
$$\lim_{n\to +\infty}\frac{|(I_{i,n})_{C_n}\backslash \Pi_n/C_n|}{|\Pi_n/C_n|}-\frac{1}{|(I_{i,n})_{C_n}|}=0.$$ 
The second term is 
$$\frac{1}{|\Pi_n/C_n|}\bigg (\sum_{j\geq 1}\frac{(e_{i,n})_j}{e_{i,n}}\big (\sum_{Q|P_i}\frac{(e^n_Q)_j}{e^n_Q}\frac{e_{i,n}}{(e_{i,n})_j}\big )-\frac{1}{e_{i,n}}\sum_{Q|P_i}\frac{e_{i,n}}{e^n_Q}\bigg )=\sum_{j\geq 1}\frac{(e_{i,n})_j}{e_{i,n}}\frac{|((I_{i,n})_{C_n})_j\backslash \Pi_n/C_n|}{|\Pi_n/C_n|}-\frac{1}{|(I_{i,n})_{C_n}|}.$$
\subsection{Stabilization of the wild inertia}\label{stabilization}
Assume from now on that $j\geq 1$. We compute $(e_{i,n})_j$ using the diagram
\begin{center}
\begin{tikzcd}
X_{\widetilde \Pi_C(n)}\arrow[bend left, dash]{rr}{\widetilde \Pi_C(1)/\widetilde \Pi_C(n)}\arrow{rr}\arrow{dr}\arrow[bend right, dash, swap]{dr}{(\Pi_n)_{C_n}} & & X_{\widetilde \Pi_C(1)}\arrow{dl}\arrow[bend left, dash]{dl}{(\Pi_1)_{C_1}}\\
& X
\end{tikzcd}
\end{center}
Write $((I_{i,n})_{C_n})_+$ for the wild inertia subgroup of $(I_{i,n})_{C_n}$ and
$$(I_{i,n})_{C_n}(1):=ker((I_{i,n})_{C_n}\rightarrow (I_{i,1})_{C_1}), \quad e_{i,n}(1):=|(I_{i,n})_{C_n}(1)|.$$
Consider the commutative diagram with exact rows
\begin{center}
\begin{tikzcd}[column sep=scriptsize, row sep=scriptsize]
0\arrow{r}&(I_{i,n})_{C_n}(1)\arrow[hook]{d}\arrow{r}& (I_{i,n})_{C_n}\arrow{r}\arrow[hook]{d}& (I_{i,1})_{C_1}\arrow{r}\arrow[hook]{d}& 0\\
0\arrow{r}&\widetilde \Pi_C(1)/\widetilde \Pi_C(n)\arrow{r}&(\Pi_n)_{C_n}\arrow{r}&(\Pi_1)_{C_1}\arrow{r}& 0.
\end{tikzcd}
\end{center}
Since $(I_{i,n})_{C_n}(1)\subseteq \widetilde \Pi_C(1)/\widetilde \Pi_C(n)$ are $\ell$-groups by Section \ref{preliminaryyyyy}(1) and $((I_{i,n})_{C_n})_j\subseteq ((I_{i,n})_{C_n})_{+}$ are p-groups by definition, we see that
\begin{enumerate}
\item the map $(I_{i,n})_{C_n}\twoheadrightarrow (I_{i,1})_{C_1}$ induces an isomorphism $$\phi_{i,n}:((I_{i,n})_{C_n})_+\simeq ((I_{i,1})_{C_1})_+;$$
\item $((I_{i,n})_{C_n})_j\cap \widetilde \Pi_C(1)/\widetilde \Pi_C(n)=1$, so that Fact \ref{lemmaserre} below yields
$$\phi_{i,n}(((I_{i,n})_{C_n})_j)=((I_{i,1})_{C_1})_{\left\lceil j/e_{i,n}(1) \right\rceil}.$$
\end{enumerate}
Write $j_{i,0}$ for smallest integer $\geq 0$ such that $(e_{i,1})_{j_{i,0}}=0$ and
$$((\widetilde I_{i,n})_{C_n})_j:=\phi_{i,n}^{-1}((I_{i,1})_j)\subseteq (I_{i,n})_{C_n}, \quad ((\widetilde I_{i})_C)_j:=\varprojlim_i ((\widetilde I_{i,n})_{C_n})_j\subseteq \Pi_C.$$
Combining (1) and (2), we get
$$((I_{i,n})_{C_n})_j=((\widetilde I_{i,n})_{C_n})_{\left\lceil j/e_{i,n}(1) \right\rceil}=$$
$$\begin{cases} ((\widetilde I_{i,n})_{C_n})_k=( (\widetilde I_i)_C)_k & \mbox{ if } \exists \mbox{ } 1\leq k\leq j_{i,0} \mbox{ such that } e_{i,n}(1)(k-1)< j \leq e_{i,n}(1)k \\ 1 & \mbox{ if } j>e_{i,n}(1)j_{i,0}.
\end{cases}.$$
\begin{fact}[{\cite[Lemma 5, IV, Pag. 75]{Serre}}]\label{lemmaserre}
Let $K\subseteq L$ a finite Galois extension of local fields with group $G$. For $-1\leq u\in \mathbb R$, write $G_u$ for the $
\left\lceil u \right\rceil ^{th}$ ramification group in lower numbering and consider the function:
$$\psi_{L/K}(u)=\int_{0}^{u} \frac{1}{[G_0:G_u]} \, \mathrm{d}t.$$
If $N\subseteq G$ if a normal subgroup corresponding to a Galois extension $K\subseteq K'$, then 
$$G_uN/N=(G/N)_{\psi_{L/K'}(u)}$$
\end{fact}
\subsection{End of proof}\label{endend}
We can continue the computation
$$\sum_{j\geq 1}\frac{|((I_{i,n})_{C_n})_j|}{|(I_{i,n})_{C_n}|}\frac{|((I_{i,n})_{C_n})_j\backslash \Pi_n/C_n|}{|\Pi_n/C_n|}-\frac{1}{|(I_{i,n})_{C_n}|}=$$
$$e_{i,n}(1)\sum_{1\leq k\leq j_{i,0}}\frac{|((\widetilde I_{i,n})_{C_n})_k||((\widetilde I_{i,n})_{C_n})_k\backslash \Pi_n/C_n|}{|(I_{i,n})_{C_n}||\Pi_n/C_n|}-\frac{1}{|(I_{i,n})_{C_n}|}=$$
$$\frac{e_{i,n}(1)}{|(I_{i,n})_{C_n}|}\sum_{1\leq k\leq j_{i,0}}
|((\widetilde I_i)_C)_k|\big (\frac{|((\widetilde I_{i,n})_{C_n})_k\backslash \Pi_n/C_n|}{|\Pi_n/C_n|}-\frac{1}{|((\widetilde I_i)_C)_k|}\big )=$$
$$\frac{1}{|(I_{i,1})_{C_1}|}\sum_{1\leq k\leq j_{i,0}}
|((\widetilde I_i)_C)_k|\big (\frac{|((\widetilde I_{i,n})_{C_n})_k\backslash \Pi_n/C_n|}{|\Pi_n/C_n|}-\frac{1}{|( (\widetilde I_i)_C)_k|}\big ).$$
Setting $(\widetilde I_{i})_k$ for the preimage of $((\widetilde I_{i})_C)_k$ under the map $\Pi\twoheadrightarrow \Pi_C$ and observing that $((\widetilde I_{i})_k)_C=((\widetilde I_{i})_C)_k$, we conclude the proof since
$$\lim_{n\to+\infty}\frac{|((\widetilde I_{i,n})_{C_n})_k\backslash \Pi_n/C_n|}{|\Pi_n/C_n|}-\frac{1}{|( (\widetilde I_i)_C)_k|}=\lim_{n\to+\infty}\frac{|((\widetilde I_{i,n})_{C_n})_k\backslash \Pi_n/C_n|}{|\Pi_n/C_n|}-\frac{1}{|( (\widetilde I_i)_k)_C|}=0$$
by Fact \ref{limite} and Section \ref{preliminaryyyyy}(1).
\section{Proof of Theorem \texorpdfstring{ \ref{UOI} }  -}\label{secmain}
\subsection{Projective systems of abstract modular scheme}
\subsubsection{Group theory}
Fix a closed subgroup $\Pi$ of $GL_r(\Zl)$, write $\Phi(\Pi)$ for the Frattini subgroup of $\Pi$, i.e. the intersection of the maximal open subgroups of $\Pi$. Set $\mathcal C_{0}(\Pi):=\{\Pi\}$ and for every integer $n\geq 1$ define $\mathcal C_{n}(\Pi)$ as the set of open subgroups $U\subseteq \Pi$ such that $\Phi(\Pi(n-1)) \subseteq U$ but $\Pi(n-1) \not\subseteq U$. By \cite[Lemma 3.1]{UOI1}, the maps $\psi_{n}:\calC_{n+1}(\Pi)\rightarrow \calC_{n}(\Pi)$ $\psi_{n}:U\mapsto U\Phi(\Pi(n-1))$ are well defined and they endow to the collection $\{\mathcal C_{n}(\Pi)\}_{n\in \mathbb N}$ with a structure of a projective system. For any $\underline C:=(C[n])_{n\geq 0}\in\varprojlim \calC_{n}(\Pi)$ write $$C[\infty]:=\varprojlim C[n]=\cap C[n]\subseteq \Pi.$$ By \cite[Lemma 3.3]{UOI1}, one has the following.
\begin{lemma}
\begin{enumerate}\label{propertyprojectivesistem}
\item[]
\item $\calC_{n}(\Pi)$ is finite and, for $n\gg 0$ (depending only on $\Pi$), it coincide with set of open subgroups $U\subseteq \Pi$ such that $\Pi(n) \subseteq U$ but $\Pi(n-1) \not\subseteq U$
\item For any $\underline C:=(C[n])_{n\geq 0}\in \varprojlim_n \mathcal C_{n}(\Pi)$, the subgroup $C[\infty]$ is a closed but not open subgroup of $\Pi$.
\item For any closed subgroup $C\subseteq \Pi$ such that $\Pi(n-1)\not\subseteq C$ there exists $U\in \mathcal C_{n}(\Pi)$ such that $C\subseteq U$.
\end{enumerate}
\end{lemma}
\subsubsection{Anabelian dictionary}\label{definitionprojective}
Let $X$ be a smooth geometrically connected $k$-variety and assume now that $\Pi$ is the image of a continuous representation $\rho:\pi_1(X)\rightarrow GL_r(\Zl)$. Consider the following (possibly disconnected) \etale covers:
$$f_{n}:\mathcal X_{n}:=\coprod_{U\in \calC_{n}(\Pi)} X_U\rightarrow X.$$
\begin{proposition}\label{lemmafondamentale}
Let $n$ be an integer $\gg 0$ (depending only on $\Pi$). If $x\in X(k)-f_{n}(\mathcal X_{n}(k))$, then $\Pi(n-1)\subseteq \Pi_x$, hence $[\Pi:\Pi_x]\leq [\Pi:\Pi(n-1)]$.
\end{proposition}
\proof
This follows from Fact \ref{formalprop} and Lemma \ref{propertyprojectivesistem}(3).
\endproof
Assume from now on that $X$ is a curve. From Theorem \ref{genus} we deduce:
\begin{corollary}\label{gonalityinteressante}
Assume that $\rho$ is GLP, $\ell\neq p$ and fix two integers $d_1,d_2\geq 1$. Then there exists an integer $N\geq 1$, depending only on $\rho,d_1,d_2$, such that for every $n\geq N$ and every $U\in \calC_{n}(\Pi)$ we have $[k_U:k]>d_1$ or $g_U>d_2$. 
\end{corollary}
\proof
This follows from Theorem \ref{genus} arguing as in \cite[Corollaries 3.7 and 3.8]{UOI1}.
\endproof
\subsection{Proof of Theorem \texorpdfstring{ \ref{UOI} }  - and a corollary}
\subsubsection{Proof of Theorem \texorpdfstring{ \ref{UOI} }  -}
Assume that $X$ is a curve and $\rho$ is $GLP$.
Consider the projective system of covers constructed in \ref{definitionprojective}
$$f_{n}:\mathcal X_{n}:=\coprod_{U\in \calC_{n}(\Pi)} X_U\rightarrow X.$$
By Corollary \ref{gonalityinteressante} we can choose an $n_0$ such that each connected component of $\mathcal X_{n_0}$ has genus lager then the constant $g$ of Fact \ref{mordel} or is defined over a non trivial extension of $k$. By the choice of $n_0$, the image $X_{n_0}$ of $f_{n_0}:\mathcal X_{n_0}(k)\rightarrow X(k)$ has a finite number of points. Replacing $n_0$ with some integer $n_0'\geq n_0$, by Lemma \ref{lemmafondamentale} for all $x\in X(k){-}X_{n_0}$ we have $\Pi(n_0)\subseteq \Pi_x$. Hence $X^{ex}(k)\subseteq X_{n_0}$ is finite and one can take 
$$N := \max_{x\in X_{n_0}{-}X^{ex}_{\rho}(k)}{\Big \{[\Pi : \Pi(n_0)], [\Pi : \Pi_x]\Big \}}.$$
This concludes the proof of Theorem \ref{UOI}.
\subsubsection{Uniform boundedness of \texorpdfstring{$\ell$-}-primary torsion}
For further use, we state a generalization of Corollary \ref{uniformintroab} for arbitrary GLP representations. We recall the notation and the terminology of \cite[Section 4]{UOI1}. Given a finitely generated free $\Zl$ module $T\simeq \Zl^r$ with a continuous action of $\pi_1(X)$ write $V:=T\otimes \Q_{\ell}$ and $M:=V/T$. For a character $\chi:\pi_1(X)\rightarrow \Zl^*$, a field extension $k\subseteq L$ and a morphism $\xi: Spec(L)\rightarrow X$, let $\chi_\xi$ (resp. $\rho_{\xi}$) denote the composition of $\chi$ (resp. $\rho$) with the morphism $\pi_1(L)\rightarrow \pi_1(X)$. Consider the following $\pi_1(L)$-sets
$$\overline M_{\xi}:=\{v\in M \ | \ \rho_\xi(\sigma)v\in <v>\},\quad \overline T_{\xi}:=\{v\in T \ | \ \rho_\xi(\sigma)v\in <v>\},$$
and $\pi_1(L)$-modules
$$M_{\xi}(\chi):=\{v\in M \ | \ \rho_\xi(\sigma)v=\chi_\xi(\sigma)v\}, \quad T_{\xi}(\chi):=\{v\in T \ | \ \rho_\xi(\sigma)v=\chi_\xi(\sigma)v\}.$$
Recall that $\chi$ is said to be non-sub-$\rho$ if $\chi_{x}$ is not isomorphic to a sub representation of $\rho_x$ for any $x\in X(k)$.
Finally denote with $T_{(0)}$ the maximal isotrivial submodule of $T$, i.e. the maximal submodule of $T$ on which
$\pi_1(X_{\overline k})$ acts via a finite quotient.
\begin{corollary}\label{boundedness}
Assume that $k$ is  finitely generated, $X$ is a curve, $\ell\neq p$ and that $\rho:\pi_1(X)\rightarrow GL(T)$ is $GLP$. Then
\begin{enumerate}
\item For every non-sub-$\rho$ character $\chi:\pi_1(X)\rightarrow \Zl^*$, there exists an integer $N\geq 1$, depending only on $\rho$ and $\chi$, such that, for any $x\in X(k)$ the $\pi_1(x)$-module $M_x(\chi)$ is contained in $M[\ell^N]$.
\item Assume furthermore that $T_{(0)}=0$. Then there exists an integer $N\geq 1$, depending only on $\rho$, such that for every $x\in X(k)-X^{ex}_{\rho}(k)$, the $\pi_1(k)$-set $\overline M_x$ is contained in $M[\ell^N]$.
\end{enumerate}
\end{corollary}
\proof
This follows from Theorem \ref{UOI} as in the proof of \cite[Corollary 4.3]{UOI1}.\endproof
\subsection{Further remarks}\label{secquestion1}
Let $k$ be a finitely generated field of characteristic $p\geq 0$, let $X$ be a smooth geometrically connected $k$-curve.
Let $\rho:\pi_1(X)\rightarrow GL_r(\Zl)$ be a continuous representation and retain the notation of Section \ref{inex}.
\subsubsection{Points of bounded degree}
As already mentioned in Section \ref{statementUOI}, Theorem \ref{UOI} is the natural extension to positive characteristic of the main result of \cite{UOI1}. In the subsequent paper \cite{UOI2}, Cadoret-Tamagawa show (\cite[Theorem 1.1]{UOI2}) that if $p=0$ and $\rho$ is GLP, then for every $d\geq 1$, the set $X_{\rho}^{ex}(\leq d)$ is finite and there exists an integer $N(\rho,d):=N\geq 1$, depending only on $\rho$ and $d$, such that $[\Pi:\Pi_x]\leq N$ for all $x\in X_{\rho}^{gen}(\leq d)$. To prove this they study the gonality of the connected components of the abstract modular curves $\mathcal X_n$.
\subsubsection{Gonality}
For a smooth proper $k$-curve $Y$, the (geometric) gonality $\gamma_Y$ of $Y$ is the minimum degree of a non constant map $Y_{\overline k}\rightarrow \mathbb P_{\overline k}^1$. While the genus $g_Y$ controls the finiteness of $k$-rational points, the gonality, in characteristic zero, controls the finiteness of points of bounded degree.
\begin{fact}[\cite{faltings}, \cite{frey}]\label{isogonality}
If $k$ is a finitely generated field of characteristic zeros and $d\geq 1$ is an integer, for every smooth proper $k$-curve $Y$ such that $\gamma_Y \geq 2d + 1$, the set $Y(\leq d)$ is finite.
\end{fact}
In light of Fact \ref{isogonality} and of the strategy described in Section \ref{ingenusgonality}, to prove \cite[Theorem 1.1]{UOI2} when $p=0$, Cadoret-Tamagawa show (\cite[Theorem 3.3]{UOI2}) that, for every $C\subseteq \Pi$ closed but not open subgroup, the gonality $\gamma_{C\Pi_{\overline k}(n)}$ of the smooth compactification of $X_{C\Pi_{\overline k}(n)}$ tends to infinity with $n$. While one can adapt (\cite{miotesi}) the arguments of  \cite[Theorem 3.3]{UOI2} to prove that $\gamma_{C\Pi_{\overline k}(n)}$ tends to infinity also when $p>0$, the positive characteristic variant\footnote{This is due to isotriviality issues in the positive characteristic variant of the Mordell-Lang conjecture; see \cite[Appendix]{JN}.} of Fact \ref{isogonality} is not true, so that one cannot deduces directly from the growth of the gonality the positive characteristic analogue of \cite[Theorem 1.1]{UOI2}.
\subsubsection{Isogonality}
However, in \cite[Appendix]{JN} Cadoret-Tamagawa have introduced a new invariant, the isogonality, that could be used to study points of bounded degree is positive characteristic. 
\begin{definition}
Let $k$ a field of characteristic $p > 0$ and $Y$ a smooth proper geometrically connected $k$-curve. The $\overline k$-isogonality $\gamma^{iso}_{Y}$ of $Y$ is defined as $d + 1$, where $d$ is the smallest integer which satisfies the following condition:
\begin{itemize}
\item There is no diagram $Y_{\overline k}\leftarrow Y'\rightarrow B$
of non constant morphisms of smooth proper curves over $\overline k$, with $B$ an
isotrivial\footnote{If $k$ is a field of characteristic $p>0$, a $k$-scheme $S$ is said to be isotrivial, if there exists a finite field $\F_q \subseteq k$ and a $\F_q$-scheme $S_0$ such that $S_{0}\times_{\F_q}\overline k\simeq S_{\overline k}$.} curve and $deg(Y'\rightarrow B)\leq d$.
\end{itemize}
\end{definition}
Their result is the following:
\begin{fact}[{\cite[Corollary A.7]{JN}}]
If $k$ is a finitely generated field of positive characteristic and $d\geq 1$ is an integer,  and $d\geq 1$ is an integer, then for every smooth proper $k$-curve $Y$ such that $\gamma_Y \geq 2d + 1$ and $\gamma^{iso}_Y\geq d+1$, the set $Y(\leq d)$ is finite.
\end{fact}
Since, by \cite{miotesi}, we know that $\gamma_{C\Pi_{\overline k}(n)}$ tends to infinity, to extend Theorem \ref{UOI} to points of bounded degree it would be enough to show the following.
\begin{conjecture}
Assume that $\rho$ is a GLP, $p>0$ and $\ell\neq p$. Then for every closed but not open subgroup
$C\subseteq \Pi_{\overline k}$ one has
$$\lim_{n\to +\infty}\gamma^{iso}_{X_{C\Pi_{\overline k}(n)}}=+\infty.$$
\end{conjecture}
 
\vskip 1\baselineskip
 \textit{emilianomaria.ambrosi@gmail.com}\\
Centre de Math\'ematiques Laurent Schwartz - Ecole Polytechnique,\\
91128 PALAISEAU, FRANCE.
\end{document}